\documentclass{amsart}
\usepackage{amssymb,amsfonts,amsbsy,url}
\usepackage[hypertex]{hyperref}
\hypersetup{colorlinks,%
citecolor=blue,%
linkcolor=blue,%
urlcolor=blue,%
}

\def \C {\mathbb{C}}
\def \Q {\mathbb{Q}}
\def \ni {\noindent}

\def \ms {\medskip}
\def \mni {\ms\ni}

\def \A {\mathcal{A}}
\def \X {\mathcal{X}}
\def \Y {\mathcal{Y}}

\def \S {\mathcal{S}}
\def \T {\mathcal{T}}

\def \K {\mathcal{K}}
\def \P {\mathcal{P}}
\def \ci {\mathcal{I}}

\def \exp {{\rm exp}}
\def \deg {{\rm deg}}

\def \jac {{\rm Jac}}

\def\Maple{{\sffamily\small Maple\/ }}

\newtheorem{theorem}{Theorem}[section]
\newtheorem{lemma}[theorem]{Lemma}

\newtheorem{definition}[theorem]{Definition}
\newtheorem{example}[theorem]{Example}

\newtheorem{corollary}[theorem]{Corollary}

\newtheorem{remark}[theorem]{Remark}

\numberwithin{equation}{section}

\begin{document}

\title[Locally nilpotent derivations in dimension three]
{Triangulable locally nilpotent derivations in dimension three}

\author[M. A. Barkatou]{Moulay A. Barkatou}
\address{Laboratoire XLIM, UMR 6172 CNRS-Universit\'e de
Limoges\\
Avenue Albert-Thomas 123, 87060, Limoges Cedex, France}
\email{moulay.barkatou@unilim.fr}

\author[H. El Houari]{Hassan El Houari}
\address{Department of Mathematics\\Faculty of Sciences Semlalia\\
Cadi Ayyad University\\P.O Box 2390, Marrakech\\Morocco}
\email{h.elhouari@ucam.ac.ma}

\author[M. El Kahoui]{M'hammed El Kahoui}
\address{Department of Mathematics\\Faculty of Sciences Semlalia\\
Cadi Ayyad University\\P.O Box 2390, Marrakech\\Morocco}
\email{elkahoui@ucam.ac.ma}

\date{}

\dedicatory{}



\begin{abstract}In this paper we give an algorithm to
recognize triangulable  locally nilpotent derivations in dimension
three. In case the given derivation is triangulable, our method
produces a coordinate system in which it exhibits a triangular
form.
\end{abstract}

\maketitle


\section{Introduction}\label{sec:intro}
Let $\K$ be a commutative field of characteristic zero, $\K^{[n]}$
be the ring of polynomials in $n$ variables with coefficients in
$\K$ and $Aut_{\K}(\K^{[n]})$ be the group of $\K$-automorphisms
of $\K^{[n]}$. Let $\underline{x}=x_1,\ldots,x_n$ be a coordinate
system of $\K^{[n]}$, i.e., $\K^{[n]}=\K[x_1,\ldots,x_n]$. Then
any automorphism $\sigma \in Aut_{\K}(\K^{[n]})$ is uniquely
determined by the images $\sigma(x_1),\ldots,\sigma(x_n)$. The
{\it affine subgroup} of $Aut_{\K}(\K^{[n]},\underline{x})$ with
respect to $\underline{x}$ is defined as
$$Af_{\K}(\K^{[n]},\underline{x})=
\{\sigma\; ; \; \deg(\sigma(x_i))=1,\; i=1,\ldots,n\}.$$ The {\it
triangular subgroup} of $Aut_{\K}(\K^{[n]})$ with respect to
$\underline{x}$ is defined as
$$BA_{\K}(\K^{[n]},\underline{x})=\{\sigma\; ; \;
\sigma(x_i)=a_ix_i+f(x_1,\ldots,x_{i-1})\;, a_i\in\K^{\star},\;
i=1,\ldots,n\}.$$ The {\it tame subgroup} of $Aut_{\K}(\K^{[n]})$
with respect to $\underline{x}$ is the subgroup generated by
affine and triangular automorphisms, and is denoted by
$TA_{\K}(\K^{[n]},\underline{x})$. Automorphisms which belong to
$TA_{\K}(\K^{[n]},\underline{x})$ are called {\it tame}, and those
which are not tame are called {\it wild}.

\medskip Automorphisms of $\K^{[2]}$ are well understood. They are
all tame and $Aut_{\K}(\K^{[2]})$ is the free amalgamated product
of $Af_{\K}(\K^{[2]},x_1,x_2)$ and $BA_{\K}(\K^{[2]},x_1,x_2)$
along their intersection \cite{jung42,kulk53}. But so far
$Aut_{\K}(\K^{[n]})$ remains a big mystery for $n\geq 3$, and it
is only recently that the existence of wild automorphisms was
established \cite{shestakov2004a}.

\medskip In order to understand the nature of $Aut_{\K}(\K^{[n]})$
it is natural to investigate algebraic group actions on the affine
$n$-space over $\K$. Actions of the algebraic group $(\K,+)$ are
commonly called {\it algebraic $G_a$-actions}, and are of the form
$\exp(t\X)_{t\in \K}$ where $\X$ is a locally nilpotent
$\K$-derivation of the polynomial ring $\K^{[n]}$.

\medskip A locally nilpotent $\K$-derivation $\X$ of
$\K[\underline{x}]$ is called {\it triangular} in the coordinate
system $\underline{x}$ if for any $i=1,\ldots,n$ we have
$\X(x_i)\in \K[x_1,\ldots,x_{i-1}]$. This is equivalent to the
fact that its one-parameter group $\exp(t\X)_{t\in \K}$ is a
subgroup of $BA_{\K}(\K^{[n]},\underline{x})$. The $\K$-derivation
$\X$ is called {\it triangulable} if there exists a
$\K$-automorphism $\sigma$ of $\K^{[n]}$ such that
$\sigma\X\sigma^{-1}$ is triangular in the coordinate system
$\underline{x}$, i.e., there exists a coordinate system
$\underline{u}$ in which $\X$ has a triangular form.

\medskip A natural question is to decide whether a given locally
nilpotent $\K$-derivation is triangulable. Bass was the first to
give in \cite{bass84a} an example of non-triangulable locally
nilpotent derivation in dimension $3$. Bass' construction was
generalized by Popov in \cite{popov87} to obtain non-triangulable
locally nilpotent derivations in any dimension $n\geq 3$. A
necessary condition of triangulability, based on the structure of
the variety of fixed points, is also given in \cite{popov87}.
However, this condition is not sufficient as proven in
\cite{daigle96}. Other criteria of triangulability in dimension
$3$ are given in
\cite{freudenburg95a,daigle98a,daigle96,walcher97a}. But it is
nowhere near obvious to make them working in an algorithmic way.

\medskip The aim of the present paper is to develop
an algorithm to check whether a given locally nilpotent derivation
$\X$ of $\K[x,y,z]$ is triangulable, and if so to find a
coordinate system $u,v,w$ in which $\X$ has a triangular form.

\medskip The paper is structured as follows. In section
\ref{sec:basics} we recall the basic facts on locally nilpotent
derivations and coordinates to be used in the paper. In section
\ref{sec:characterization} we give an algorithmic characterization
of rank two locally nilpotent derivations in dimension three. A
triangulability criterion is given in section
\ref{sec:triangulability}, while section \ref{sec:triang_system}
contains the main ingredients that make this criterion work in an
algorithmic way. Computational examples top off the paper.

\section{Notation and basic facts}\label{sec:basics}
Throughout this paper $\K$ is a commutative field of
characteristic zero, all the considered rings are commutative of
characteristic zero with unit and all the considered derivations
are nonzero. A derivation of a $\K$-algebra $\A$ is called a
$\K$-derivation if it satisfies $\X(a)=0$ for any $a\in \K$.

\subsection{Coordinates}
A polynomial $f\in \K[x_1,\ldots,x_n]$ is called a {\it
coordinate} if there exists a list of polynomials $f_1,\ldots,
f_{n-1}$ such that $\K[x_1,\ldots, x_n]=\K[f,f_1,\ldots,
f_{n-1}]$. A list $f_1,\ldots, f_r$ of polynomials, with $r\leq
n$, is called a {\it system of coordinates} if there exists a list
$f_{r+1},\ldots,f_n$ of polynomials such that
$\K[x_1,\ldots,x_n]=\K[f_1,\ldots,f_n]$. A system of coordinates
of length $n$ will be called a {\it coordinate system}.

\mni The Abhyankar-Moh Theorem \cite{abhyankar-moh75a} states that
a polynomial $f$ in $\K[x,y]$ is a coordinate if and only if
$\K[x,y]/f$ is $\K$-isomorphic to $\K^{[1]}$. In the case of three
variables we have the following result proved by Kaliman in
\cite{kaliman2002a} for the case $\K=\C$ and extended to the case
of arbitrary commutative fields of characteristic zero in
\cite{daigle-kaliman2004a}.
\begin{theorem}\label{kaliman-daigle}Let $f$ be a polynomial in
$\K[x,y,z]$ and assume that for all but finitely many $\alpha\in
\K$ the $\K$-algebra $\K[x,y,z]/(f-\alpha)$ is $\K$-isomorphic to
$\K^{[2]}$. Then $f$ is a coordinate of $\K[x,y,z]$.
\end{theorem}
\mni A given polynomial $f$ of $\K[x_1,\ldots, x_n]$ is called a
{\it local coordinate} if it satisfies
$\K(f)[x_1,\ldots,x_n]\simeq_{\K(f)}\K(f)^{[n-1]}$. As a
consequence of Theorem \ref{kaliman-daigle}, any local coordinate
of $\K[x,y,z]$ is in fact a coordinate, see \cite{elkahoui2005a}.
The original proof of Theorem \ref{kaliman-daigle} is of
topological nature, and it is not clear how to compute polynomials
$g,h$ such that $\K[f,g,h]=\K[x,y,z]$.

\medskip The study of coordinates in polynomial rings over fields
naturally leads to do the same but over rings. Given a ring $\A$
and $f\in \A[x_1,\ldots,x_n]$, we say that $f$ is a residual
coordinate if $f$ is a coordinate of $\K_{\P}[x_1,\ldots,x_n]$ for
any prime ideal $\P$ of $\A$, where $\K_{\P}$ stands for the
residual field of $\A$ in $\P$. The following result, proved in
\cite{bhat-duta93a} for the Noetherian case and extended to the
general case in \cite{essen2004a}, will be crucial for our
purpose.

\begin{theorem}\label{residual_coordinate}Let $\A$ be a ring
containing $\Q$. Then any residual coordinate of $\A[x,y]$ is a
coordinate of $\A[x,y]$.
\end{theorem}

\subsection{Locally nilpotent derivations}
A derivation of a ring $\A$ is called {\it locally nilpotent} if
for any $a\in\A$ there exists a positive integer $n$ such that
$\X^n(a)=0$. The subset $\{a\in \A\; ; \;\X(a)=0\}$ of $\A$ is in
fact a subring called the {\it ring of constants} of $\X$ and is
denoted by $\A^{\X}$. When $\A$ is a domain and $\X$ is locally
nilpotent, the ring of constants $\A^{\X}$ is factorially closed
in $\A$, i.e., if $a\in \A^{\X}$ and $a=bc$ then $b,c\in \A^{\X}$.
In particular the units of $\A$ are in $\A^{\X}$ and the
irreducible elements of $\A^{\X}$ are irreducible in $\A$.

\medskip An element $s$ of $\A$ satisfying $\X(s)\neq 0$ and
$\X^2(s)=0$ is called a {\it local slice} of $\X$. If moreover
$\X(s)=1$ then $s$ is called a {\it slice} of $\X$. A locally
nilpotent derivation needs not to have a slice but always has a
local slice. The following result, which dates back at least to
\cite{wright81a}, concerns locally nilpotent derivations having a
slice.
\begin{lemma}\label{wright}Let $\A$ be a ring containing $\Q$
and $\X$ be a locally nilpotent derivation of $\A$ having a slice
$s$. Then $\A=\A^{\X}[s]$ and $\X={\partial}_s$.
\end{lemma}

\mni Locally nilpotent derivations in two variables over fields
are well understood. We have in particular the following version
of Rentschler's Theorem \cite{rentschler68}.

\begin{theorem}\label{rentshcler}Let $\X$ be a locally nilpotent
$\K$-derivation of $\K[x,y]$. Then there exists a coordinate
system $f,g$ of $\K[x,y]$ and a univariate polynomial $h$ such
that $\K[x,y]^{\X}=\K[f]$ and $\X=h(f){\partial}_g$.
\end{theorem}

\medskip  As a consequence of Theorem \ref{rentshcler}, if $\A$ is a
UFD containing $\Q$ and $\X$ is a locally nilpotent
$\A$-derivation of $\A[x,y]$ then there exists $f\in \A[x,y]$ and
a univariate polynomial $h$ such that $\A[x,y]^{\X}=\A[f]$ and
$\X=h(f)({\partial}_yf{\partial}_x-{\partial}_xf{\partial}_y)$,
see \cite{daigle98a}. In case $\A$ is an arbitrary ring, the
situation is much more involved, see e.g., \cite{bhatwadekar97a}.
However, we have the following result from \cite{berson2001a}.
\begin{theorem}\label{berson_theorem}Let $\A$ be a ring containing
$\Q$ and $\X$ be a locally nilpotent $\A$-derivation of $\A[x,y]$
such that $1\in \ci(\X(x),\X(y))$. Then there exists a polynomial
$f$ such that $\A[x,y]^{\X}=\A[f]$ and $\X$ has a slice $s$. In
particular, $\A[x,y]=\A[f,s]$ and $\X={\partial}_s$.
\end{theorem}

\medskip In case $\A=\K^{[3]}$ we have the following result proved
by Miyanishi \cite{Miyanishi85a} for the case $\K=\C$ and extended
to the general case in a straightforward way by using Kambayashi's
result \cite{kambayashi75a}, see also \cite{makar_limanov2005} for
an algebraic proof.
\begin{theorem}\label{miyanishi}Let $\X$ be a locally nilpotent
$\K$-derivation of $\K[x,y,z]$. Then there exist $f,g\in
\K[x,y,z]$ such that $\K[x,y,z]^{\X}=\K[f,g]$.
\end{theorem}

\subsection{Rank of a derivation}
Let $\X$ be a $\K$-derivation of $\K[\underline{x}]=\K[x_1,\ldots,
x_n]$. As defined in \cite{freudenburg95a} the {\it co-rank} of
$\X$, denoted by $corank(\X)$, is the unique nonnegative integer
$r$ such that $\K[\underline{x}]^{\X}$ contains a system of
coordinates  of length $r$ and no system of coordinates of length
greater than $r$. The {\it rank} of $\X$, denoted by $rank(\X)$,
is defined by $rank(\X)=n-corank(\X)$. Intuitively, the rank of
$\X$ is the minimal number of partial derivatives needed for
expressing $\X$. The only one derivation of rank $0$ is the zero
derivation. Any $\K$-derivation of rank $1$ is of the form
$p(f_1,\ldots,f_{n}){\partial}_{f_n}$, where $f_1,\ldots,f_n$ is a
coordinate system. Such a derivation is locally nilpotent if and
only if $p$ does not depend on $f_n$.

\medskip Let $\X$ be a locally nilpotent $\K$-derivation of
$\K[\underline{x}]$ and let us consider
$c=\gcd(\X(x_1),\ldots,\X(x_n))$. We say that $\X$ is {\it
irreducible} if $c$ is a constant of $\K^{\star}$. It is well
known that $\X(c)=0$ and $\X=c\Y$ where $\Y$ is an irreducible
locally nilpotent $\K$-derivation. Moreover, this decomposition is
unique up to a unit, i.e., if $\X=c_1\Y_1$, where $\Y_1$ is
irreducible, then there exists a constant $\mu\in \K^{\star}$ such
that $c_1=\mu c$ and $\Y=\mu\Y_1$.

\medskip Given any irreducible locally nilpotent $\K$-derivation of
$\K[x_1,\ldots, x_n]$ and any $c$ such that $\X(c)=0$, the
derivations $\X$ and $c\X$ have the same rank. Thus, for rank
computation we may reduce, without loss of generality, to
irreducible derivations. We will see in section
\ref{sec:characterization} that the rank of a locally nilpotent
derivation in dimension three may be computed by using classical
techniques of computational commutative algebra, namely Gr\"obner
bases and functional decomposition of multivariate polynomials.

\subsection{The plinth ideal and minimal local slices}Let
$\A$ be a ring, $\X$ be a locally nilpotent derivation of $\A$ and
let
$$\S^{\X}:=\{\X(a)\; ;\; \X^{2}(a)=0\}.$$ It is easy to see
that $\S^{\X}$ is an ideal of $\A^{\X}$, called the {\it plinth
ideal} of $\X$. This is clearly an invariant of $\X$, i.e.,
$\S^{\sigma\X\sigma^{-1}}=\sigma(\S^{\X})$ for any automorphism
$\sigma$ of $\A$. In case $\A=\K[x,y,z]$ we have the following
result which is a direct consequence of faithful flatness of
$\K[x,y,z]$ over $\K[x,y,z]^{\X}$, see \cite{daigle-kaliman2004a}
for the general case.
\begin{theorem}\label{bonnet}Let $\X$ be a locally nilpotent
$\K$-derivation of $\K[x,y,z]$. Then the plinth ideal $\S^{\X}$ is
principal.
\end{theorem}

\medskip Computing a generator of the ideal $\S^{\X}$ is of central
importance for our purpose. For this we need the concept of {\it
minimal local slice} which may be found in
\cite{freudenburg_book,elkahoui2006a}.

\begin{definition}Let $\A$ be a domain and $\mathcal{X}$ be a locally
nilpotent derivation of $\A$. A local slice $s$ of $\mathcal{X}$
is called minimal if for any local slice $v$ such that
$\mathcal{X}(v)\;\vert\; \mathcal{X}(s)$ we have
$\mathcal{X}(v)=\mu \mathcal{X}(s)$, where $\mu$ is a unit of
$\A$.
\end{definition}

\begin{lemma}\label{local_slice}Let $\A$ be a UFD, $\X$ be a
locally nilpotent derivation of $\A$ and $s$ be a local slice of
$\X$. Then the following hold:

i) there exists a minimal local slice $s_0$ of $\X$ such that
$\X(s_0)\;\vert\; \X(s)$,

ii) in case $\S^{\X}$ is a principal ideal, it is generated by
$\X(s)$ for any minimal local slice $s$ of $\X$.
\end{lemma}
\begin{proof}$i)$ Let $s$ be a local slice of $\X$ and write
$\X(s)=\mu p_1^{m_1}\cdots p_r^{m_r}$, where $\mu$ is a unit and
the $p_i$'s are prime, and set $m=\sum_im_i$. We will prove the
result by induction on $m$.

\mni For $m=0$ we have $\X(s)=\mu$, and so $\mu^{-1}s$ is a slice
of $\X$. This shows that $s$ is a minimal local slice of $\X$. Let
us now assume the result to hold for $m-1$ and let $s$ be a local
slice of $\X$, with $\X(s)=\mu p_1^{m_1}\cdots p_r^{m_r}$ and
$\sum_im_i=m$. Then we have one of the following cases.

\mni -- For any $i=1,\ldots,r$ the ideal $p_i\A$ does not contain
any element of the form $s+a$ with $\X(a)=0$. In this case $s$ is
a minimal local slice of $\X$. Indeed, if it is not the case there
exists a local slice $s_0$ of $\X$ such that $\X(s)=q\X(s_0)$,
where $q$ is  not a unit of $\A^{\X}$. Without loss of generality
we may assume that $p_1\;\vert \;q$. If we write $q=p_1q_1$ then
$\X(s-q_1p_1s_0)=0$ and so the ideal $p_1\A$ contains an element
of the form $s+a$ with $\X(a)=0$, and this contradicts our
assumption.

\mni -- There exists $i$ such that $p_i\A$ contains an element of
the form $s+a$, with $\X(a)=0$. Without loss of generality we may
assume that $i=1$. If we write $s+a=p_1s_1$ then $\X(s_1)=\mu
p_1^{m_1-1}p_2^{m_2}\cdots p_r^{m_r}$, and by using induction
hypothesis we get a minimal local slice $s_0$ of $\X$ such that
$\X(s_0)\;\vert\; \X(s_1)$. Since $\X(s_1)\;\vert \;\X(s)$ we get
the result in this case.

\medskip $ii)$ Assume now that $\S^{\X}$ is principal and let $c$ be
a generator of this ideal, with $c=\X(s_0)$. Let $s$ be a minimal
local slice of $\X$. Since $\X(s)\in \S^{\X}$ we may write
$\X(s)=\mu\X(s_0)$. The fact that $s$ is minimal implies that
$\mu$ is a unit of $\A^{\X}$, and so $\X(s)$ generates $\S^{\X}$.
\end{proof}

\medskip An algorithm for computing a generator of $\S^{\X}$ in
dimension three is given in \cite{elkahoui2006a}. As we will see
in section \ref{sec:characterization}, a generator of the ideal
$\S^{\X}$ contains crucial information for computing the rank of a
locally nilpotent derivation in dimension three.

\medskip Let $\A^{\X}[s\; ;\; \X(s)\in \S^{\X}]$ be the subring of
$\A$ generated over $\A^{\X}$ by all the local slices of $\X$.
This is another invariant of the derivation $\X$. Let $(c_i)_{i\in
I}$ be a generating system of $\S^{\X}$ and let $s_i$ be such that
$\X(s_i)=c_i$. Given any local slice $s$ of $\X$ we have $\X(s)\in
\S^{\X}$, and so there exist a finite subset $J$ of $I$ an a
family $(u_i)_{i\in J}$ in $\A^{\X}$ such that
$\X(s)=\sum_iu_i\X(s_i)$. We have then $\X(s-\sum_iu_is_i)=0$ and
so $s\in \A^{\X}[s_i,i\in I]$. This proves that the ring
$\A^{\X}[s\; ;\; \X(s)\in \S^{\X}]=\A^{\X}[s_i, i\in I]$. In case
$\S^{\X}$ is principal we get a univariate polynomial ring
$\A^{\X}[s]$, which we will call the {\it trivializing ring} of
$\X$ and denote by $\T^{\X}$.

\medskip Assume $\A$ to be a UFD and that $\S^{\X}$ is principal and
generated  by $c=\X(s)$. For any factor $q$ of $c$ we let
$\ci_{q}^\X=q\A\cap \T^{\X}[s]$. The ideals $\ci_{q}^\X$ are in
fact invariants of the derivation and we will see in section
\ref{sec:triang_system} that they hold the essential information
needed to decide whether $\X$ is triangulable.

\section{Characterization of rank two locally nilpotent
derivations}\label{sec:characterization}Let $\X$ be an irreducible
locally nilpotent derivation of $\K[x,y,z]$ and $c$ be a generator
of its plinth ideal $\S^{\X}$. From Lemma \ref{wright}, $\X$ is of
rank $1$ if and only if $c\in \K^{\star}$. The following Theorem
from \cite{elkahoui2006a} gives a characterization of rank two
locally nilpotent derivations in dimension three.

\begin{theorem}\label{rank_two_theorem}Let $\X$ be an irreducible
locally nilpotent derivation of $\K[x,y,z]$ and assume that
$rank(\X)\neq 1$. Let us write $\K[x,y,z]^{\X}=\K[f,g]$ and
$\S^{\X}=c\K[f,g]$. Then the following are equivalent:

i) $rank(\X)=2$,

ii) $c=\ell(u)$, where $\ell$ is a univariate polynomial and $u$
is a coordinate of $\K[f,g]$,

iii) $c=\ell(u)$, where $u$ is a coordinate of $\K[x,y,z]$.
\end{theorem}
\begin{proof}$i)\Rightarrow ii)$ Assume that $rank(\X)=2$ and
let $u,v,w$ be a coordinate system such that $\X(u)=0$. The
$\K$-derivation $\X$ is therefore a $\K[u]$-derivation of
$\K[u][v,w]$, and since $\K[u]$ is a UFD there exists $p\in
\K[x,y,z]$ such that $\K[f,g]=\K[u,p]$. This proves that $u$ is a
coordinate of $\K[f,g]$.

\mni Let us now view $\X$ as $\K(u)$-derivation of $\K(u)[v,w]$.
Since $\X$ is irreducible, and according to Theorem
\ref{rentshcler}, there exists $s=\frac{h(u,v,w)}{k(u)}$ such that
$\X(s)=1$, and so $\X(h)=k(u)$. Let $c$ be a generator of
$\S^{\X}$. Then $c\;\vert\; k(u)$, and since $\K[u]$ is
factorially closed in $\K[u,v,w]$ we have $c=\ell(u)$ for some
univariate polynomial $\ell$.

\medskip $ii)\Rightarrow iii)$ Assume that $c=\ell(u)$, where $u$
is a coordinate of $\K[f,g]$ and write $\K[f,g]=\K[u,p]$. Let $s$
be such that $\mathcal{X}(s)=c$. If we view $\X$ as
$\K(u)$-derivation of $\K(u)[x,y,z]$ then
$\K(u)[x,y,z]^{\X}=\K(u)[p]$ and $\X(c^{-1}s)=1$. By applying
Lemma \ref{wright} we get $\K(u)[x,y,z]=\K(u)[p,s]$. From Theorem
\ref{kaliman-daigle} we deduce that $u$ is a coordinate of
$\K[x,y,z]$.

\medskip $iii)\Rightarrow i)$ Since $rank(\X)\neq 1$ the polynomial
$\ell$ is nonconstant. We have $\X(c)=\ell^{\prime}(u)\X(u)=0$,
and so $\X(u)=0$. On the other hand, since $u$ is assumed to be a
coordinate of $\K[x,y,z]$ we have $rank(\X)\leq 2$. By assumption
we have $rank(\X)\neq 1$ and so $rank(\X)=2$.\end{proof}

\medskip The condition $ii)$ in Theorem \ref{rank_two_theorem} is
in fact algorithmic. Indeed, there are actually many algorithms to
check whether a given polynomial in two variables is a coordinate,
see e.g.,
\cite{abhyankar-moh75a,berson-essen2000a,elkahoui2004c,shpilrain97a}.
It is worth mentioning that from the complexity point of view the
algorithm given in \cite{shpilrain97a} is the most efficient as
reported in \cite{shpilrain2005a}. On the other hand, condition
$c=\ell(u)$ may be checked by using a special case, called {\it
uni-multivariate decomposition}, of functional decomposition of
polynomials, see e.g., \cite{gathen90}. It is important to notice
here that uni-multivariate decomposition is essentially unique.
Namely, if $c=\ell(u)=\ell_1(u_1)$, where $u$ and $u_1$, are
undecomposable, then there exist $\mu\in \K^{\star}$ and $\nu\in
\K$ such that $u_1=\mu u+\nu$. More details about the computation
of the rank of a locally nilpotent derivation in dimension three
may be found in \cite{elkahoui2006a}.

\section{A triangulability criterion}\label{sec:triangulability}
Triangulable derivations in dimension $n$ are of rank at most
$n-1$. On the other hand, a rank $1$ locally nilpotent derivation
is obviously triangulable. This shows that in dimension $3$ we
only need to deal with rank $2$ derivations.

\mni Let $\X$ be a rank $2$ locally nilpotent derivation of
$\K[x,y,z]$ such that $\X(x)=0$. Then for any coordinate system
$x_1,y_1,z_1$ such that $\X(x_1)=0$ we have $\K[x]=\K[x_1]$, see
\cite{daigle96} (this could also be easily deduced from the
uniqueness property of uni-multivariate decomposition). This
proves that if $\X$ has a triangular form in a coordinate system
$x_1,y_1,z_1$ then $x_1$ is essentially unique and may be
extracted from a generator of the plinth ideal $\S^{\X}$. Also,
this shows that if $\X$ is triangulable and $\X(a)=0$ then $a\X$
is triangulable if and only if $a\in \K[x]$.

\begin{lemma}\label{mini_slice_derivation}Let $\X$ be an
irreducible locally nilpotent $\K$-derivation of $\K[x,y,z]$ of
rank $\leq 2$, $u$ be a coordinate of $\K[x,y,z]$ such that
$\X(u)=0$, and $s$ be a minimal local slice of $\X$. Then the
$\K[u]$-derivation $\Y=\jac_{(x,y,z)}(u,s,.)$ is locally nilpotent
irreducible and $\K[x,y,z]^{\Y}=\K[u,s]$. Moreover, $\X\Y=\Y\X$.
\end{lemma}
\begin{proof}Without loss of generality, we may assume that $u=x$.
Let us write $\K[x,y,z]^{\X}=\K[x,p]$ and let $\X(s)=c(x)$. Then
$\K[x]_c[y,z]=\K[x]_c[p,s]$ according to Lemma \ref{wright}. Given
$a\in \K[x,y,z]$, we may therefore write
$a=\frac{h(x,p,s)}{c(x)^n}$. This gives
$\Y(a)=-c^{-n}({\partial}_zs{\partial}_yp-
{\partial}_ys{\partial}_zp){\partial}_ph$, and since
$\X(s)=-{\partial}_zs{\partial}_yp+
{\partial}_ys{\partial}_zp=c(x)$ we get
$\Y(a)=c(x)^{-n+1}\partial_ph$. By induction we get
$\Y^{d+1}(a)=0$, where $d=\deg_p(h)$, and this proves that $\Y$ is
locally nilpotent.

\medskip Let $g(x,y,z)=\gcd({\partial}_ys,{\partial}_zs)$. Since
$\Y(p)=-c(x)$ we have $g\;\vert \; c(x)$ and so we may write
$c(x)=g(x)c_1(x)$. We have then $s(x,y,z)=g(x)s_1(x,y,z)+a(x)$,
and this gives $\X(s_1)=c_1(x)$. Since $s$ is a minimal local
slice of $\X$ we have $c(x)\;\vert \; c_1(x)$, and so $g\in
\K^{\star}$. This shows that $\Y$ is irreducible.

\mni Let us write $\K[x,y,z]^{\Y}=\K[x,s_0]$ and $s=\ell(x,s_0)$.
Then $\Y=\ell^{\prime}(x,s_0)({\partial}_zs_0{\partial}_y-
{\partial}_ys_0{\partial}_z)$. Since $\Y$ is irreducible
$\ell^{\prime}$ is a unit, and so $s=\mu s_0+a(x)$. This proves
that $\K[x,s_0]=\K[x,s]$. The fact that $\X$ and $\Y$ commute is
clear.\end{proof}

\begin{lemma}\label{points_infinity}Let $\X$ be a rank two irreducible
triangulable $\K$-derivation of $\K[x,y,z]$ and let $u,v,w$ be a
coordinate system of $\K[x,y,z]$ such that
$$\X(u)=0,\quad \X(v)=d(x),\quad \X(w)=q(x,v).$$ Let $c(u)$ be a
generator of $\S^{\X}$. Then $d(u)=c(u)e(u)$, $\gcd(c(u),e(u))=1$
and $\ci(e(u),q(u,v))=\K[u,v]$.
\end{lemma}
\begin{proof}Since $\X$ is of rank $2$ we must have $d(u)\neq 0$,
and so $v$ is a local slice of $\X$. This proves that $c(u)\;
\vert \; d(u)$. On the other hand, let us consider
\begin{equation}\label{equation01}
p=d(u)w-q_1(u,v),\end{equation} where ${\partial}_vq_1=q$. We have
$\X={\partial}_wp{\partial}_v-{\partial}_vp{\partial}_w$, and the
fact that $\X$ is irreducible implies that
$\gcd({\partial}_vp,{\partial}_wp)=1$. This shows that
$\K[u,v,w]^{\X}=\K[u,p]$.

\mni Let us write $d(u)=c(u)e(u)$, and notice that the result
obviously holds if we have $\deg_u(e(u)=0$. Thus, we assume in the
rest of the proof that $\deg_u(e(u))>0$.

\mni Let $\alpha$ be a root of $e(u)$ in an algebraic closure
$\overline{\K}$ of $\K$ and let us prove that $q(\alpha,v)$ is a
nonzero constant. We may write
\begin{equation}\label{equation02}
v=e(u)s(u,v,w)+\ell(u,p(u,v,w)),\end{equation} where $s$ is a
minimal local slice of $\X$. By substituting $\alpha$ to $u$ in
the relation (\ref{equation02}) we get
$v=\ell(\alpha,p(\alpha,v,w))$, and by doing so for
(\ref{equation01}) we get $p(\alpha,v,w)=-q_1(\alpha,v)$. This
yields $v=\ell(\alpha,-q_1(\alpha,v))$. By comparing degrees in
both sides of this equality we get $\deg(q_1(\alpha,v))=1$. This
proves that $\deg(q(\alpha,v))=0$ and so $q(\alpha,v)$ is a
nonzero constant. By the Hilbert Nullstellensatz we have
$\ci(e(u),q(u,v))=\K[u,v]$. To prove that $\gcd(c,e)=1$ we only
need to show that $q(\alpha,v)$ is nonconstant for any root
$\alpha$ of $c(u)$.

\mni Let $a(u)$ be a prime factor of $c(u)$. First, notice that
the assumption $q(u,v)=0 \mod a(u)$ would imply that $a(u)\;\vert
\; \X(h)$ for any $h$ and contradicts the fact that $\X$ is
irreducible. Assume towards contradiction that $q(u,v)$ is a
nonzero constant modulo $a(u)$. Then $\X$ has no fixed points in
the surface $a(u)=0$. If we write $c(u)=a(u)^mc_1(u)$, with
$\gcd(c_1,a)=1$, and view $\X$ as $\K[u]_{c_1}$-derivation of
$\K[u]_{c_1}[v,w]$ then it is fixed point free and so it has a
slice $s$ according to Theorem \ref{berson_theorem}. If we write
$s=\frac{h(u,v,w)}{c_1^n}$ then $\X(h)=c_1^n$. But $c_1^n$ is not
a multiple of $c$, and this contradicts the fact that $c$ is a
generator of $\S^{\X}$.\end{proof}

\mni The following Lemma shows that it is possible to get rid of
the factor $e(u)$.
\begin{lemma}\label{extra_factor}Let $\X$ be a rank two irreducible
triangulable locally nilpotent $\K$-derivation of $\K[x,y,z]$, and
write $\K[x,y,z]^{\X}=\K[u,p]$ where $u$ is a coordinate of
$\K[x,y,z]$. Let $s$ be a minimal local slice of $\X$ and write
$\X(s)=c(u)$. Then there exist $v,w$ such that $u,v,w$ is a
coordinate system and
$$\X(u)=0,\;\X(v)=c(u),\;\X(w)=q(u,v).$$
\end{lemma}
\begin{proof}Let $u_1,v_1,w_1$ be a coordinate system
such that $\X(u_1)=0,\X(v_1)=d(u_1)$ and $\X(w_1)=q_1(u_1,v_1)$.
Without loss of generality we may assume that $u_1=u$, and
according to Lemma \ref{points_infinity} let us write
$d(u)=c(u)e(u)$ with $\gcd(c(u),e(u))=1$.

\mni Without loss of generality, we may choose
$p=c(u)e(u)w_1-Q_1(u,v_1)$, where ${\partial}_{v_1}Q_1=q_1$, and
$v_1=e(u)s+\ell_1(u,p)$. This gives the relation
\begin{equation}\label{equation1}
p=c(u)e(u)w_1-Q_1(u,e(u)s+\ell_1(u,p)).
\end{equation}
If we write $a(u)c(u)+b(u)e(u)=1$ then we get
$$Q_1(u,e(u)s+\ell_1(u,p))=Q_1(u,e(u)(s+b(u)\ell_1(u,p))
+c(u)a(u)\ell_1(u,p)),$$ and by Taylor expanding we get
\begin{equation}\label{equation2}
Q_1(u,e(u)s+\ell_1(u,p))=Q_1(u,e(u)(s+b(u)\ell_1(u,p)))+c(u)Q_2(u,p,s).
\end{equation}
Now, let $\ell(u,p)=b(u)\ell_1(u,p)$, $v=s+\ell(u,p)$,
$Q(u,v)=Q_1(u,e(u)v)$ and let $w=e(u)w_1-Q_2(u,p,s)$. According to
the relations (\ref{equation1}) and (\ref{equation2}) we have
\begin{equation}\label{equation3}
p+Q(u,v)=c(u)w.
\end{equation}
Let us consider the $\K[u]$-derivation $\Y=-\jac(u,v,.)$. By Lemma
\ref{mini_slice_derivation}, $\Y$ is locally nilpotent and
$\K[x,y,z]^{\Y}=\K[u,v]$. By the relation (\ref{equation3}) we
have $\Y(w)=1$, and from Lemma \ref{wright} we deduce that $u,v,w$
is a coordinate system of $\K[x,y,z]$. Moreover, we have
$\X(u)=0,\X(v)=c(u)$ and $\X(w)={\partial}_vQ(u,v)$.\end{proof}

\medskip We have now enough material to state the main result of
this section.
\begin{theorem}\label{triangulability_theorem}Let $\X$ be a rank
two irreducible locally nilpotent $\K$-derivation of $\K[x,y,z]$
and write $\K[x,y,z]^{\X}=\K[u,p]$ where $u$ is a coordinate of
$\K[x,y,z]$. Let $s$ be a minimal local slice of $\X$ and write
$\X(s)=c(u)$. Then the following are equivalent:

i) the derivation $\X$ is triangulable,

ii) the ideal $\ci_{c}^{\X}$ contains a polynomial of the form
$H=p+Q(u,s+\ell(u,p))$.

\mni In this case, if we let $v=s+\ell(u,p)$ and $H=c(u)w$ then
$u,v,w$ is a coordinate system of $\K[x,y,z]$ which satisfies
$$\X(u)=0,\;\X(v)=c(u),\;\X(w)={\partial}_vQ(u,v).$$
\end{theorem}
\begin{proof}$i)\Rightarrow ii)$ Let $u,v,w$ be a coordinate
system such that $\X(u)=0,\X(v)=d(u)$ and $\X(w)=q(u,v)$. By Lemma
\ref{extra_factor} we may choose our coordinate system in such a
way that $d(u)=c(u)$. In this case we have $v=s+\ell(u,p)$ and we
may choose $p=c(u)w-Q(u,v)$, where ${\partial}_vQ(u,v)=q(u,v)$. If
we let $H=p+Q(u,v)$ then clearly $H\in \ci_{c}^{\X}$.

\mni $ii)\Rightarrow i)$ Let $v=s+\ell(u,p)$ and
$\Y=-\jac(u,v,.)$, and notice that $\Y$ is locally nilpotent and
$\K[x,y,z]^{\Y}=\K[u,v]$ according to Lemma
\ref{mini_slice_derivation}. By assumption we have $H= p+Q(u,v)\in
\ci_{c}^{\X}$, so let us write $H=c(u)w$. Since $\Y(H)=\Y(p)=c(u)$
we have $\Y(w)=1$. According to Lemma \ref{wright}, $u,v,w$ is a
coordinate system of $\K[x,y,z]$, and $\X(u)=0,\X(v)=c(u)$ and
$\X(w)={\partial}_vQ(u,v)$.\end{proof}

\begin{corollary}\label{primary_reduction}Let $\X$ be a rank two
irreducible locally nilpotent $\K$-derivation of $\K[x,y,z]$ and
write $\K[x,y,z]^{\X}=\K[u,p]$ where $u$ is a coordinate of
$\K[x,y,z]$. Let $s$ be a minimal local slice of $\X$ and write
$\X(s)=c(u)=c_1^{n_1}\cdots c_r^{n_r}$, where the $c_i$'s are
irreducible and pairwise distinct. Then the following are
equivalent:

i) the derivation $\X$ is triangulable,

ii) for any $i=1,\ldots, r$ the ideal $\ci_{c_i^{n_i}}^{\X}$
contains a polynomial $H_i$ such that
$H_i=p+Q_i(u,s+\ell_i(u,p))\; \mod \; c_i^{n_i}$.
\end{corollary}
\begin{proof}$i)\Rightarrow ii)$ This is an obvious consequence of
Theorem \ref{triangulability_theorem}.

\medskip $ii)\Rightarrow i)$ By the Chinese remainder Theorem let
$Q(u,v)$ and $\ell(u,p)$ be such that $Q=Q_i$ and $\ell=\ell_i
\mod \; c_i^{n_i}$. A straightforward computation shows that
$p+Q(u,s+\ell(u,p))=c(u)w$, and so $\X$ is triangulable by Theorem
\ref{triangulability_theorem}.\end{proof}

\section{Computing a triangulating coordinate system}
\label{sec:triang_system}Let $\X$ be an irreducible triangulable
locally nilpotent derivation of $\K[x,y,z]$ and write
$\K[x,y,z]^{\X}=\K[u,p]$, where $u$ is a coordinate of
$\K[x,y,z]$. Let $s$ be a minimal local slice of $\X$, with
$\X(s)=c(u)=c_1^{n_1}\cdots c_r^{n_r}$ and the $c_i$'s are prime
and pairwise distinct. According to Corollary
\ref{primary_reduction} it suffices to find a polynomial of the
form $p+Q_i(u,\ell_i(u,p))$ in each ideal $\ci_{c_i^{n_i}}^{\X}$.
It is trivial to see that such a polynomial is a coordinate of
$\K[u,p,s]$, and as a by-product it is a coordinate when viewed as
polynomial of $\K[u]/c_i^{n_i}[p,s]$. We are thus led to deal with
the problem of finding a polynomial in $\ci_{c_i^{n_i}}$ which is
a coordinate of $\K[x]/c_i^{n_i}[p,s]$. In fact, taking into
account Theorem \ref{residual_coordinate}, we only need to deal
with the case of $\K[x]/c_i[p,s]$. In this section we solve such a
problem, and we show how this allows to compute a coordinate
system of $\K[x,y,z]$ in which $\X$ exhibits a triangular form.

\begin{lemma}\label{monic_lemma}Let $\X$ be a rank two
irreducible locally nilpotent $\K$-derivation of $\K[x,y,z]$, and
write $\K[x,y,z]^{\X}=\K[u,p]$ where $u$ is a coordinate of
$\K[x,y,z]$. Let $s$ be a minimal local slice of $\X$ and write
$\X(s)=c(u)$. Then for any prime factor $c_1$ of $c$ the following
hold:

i) there exists a monic polynomial $h_1$ with respect to $s$ such
that $\ci_{c_1}^{\X}=(c_1,h_1)$. Moreover, $c_1,h_1$ is the
reduced Gr\"obner basis of $\ci_{c_1}^{\X}$ with respect to the
lex-order $u\prec p\prec s$,

ii) the ideal $\ci_{c_1}^{\X}$ contains a coordinate of
$\K[u]/c_1[p,s]$ if and only if $h_1$ is a coordinate of
$\K[u]/c_1[p,s]$. Moreover, any polynomial $h\in\ci_{c_1}^{\X}$
which is a coordinate of $\K[u]/c_1[p,s]$ satisfies $h=\mu(u)h_1$,
where $\mu$ is a unit of $\K[u]/c_1$.
\end{lemma}
\begin{proof}$i)$ Let $v,w$ be such that $u,v,w$ is a
coordinate system of $\K[x,y,z]$. The derivation $\X$ induces a
locally nilpotent $\K[u]/c_1$-derivation $\overline{\X}$ of
$\K[x,y,z]/c_1=\K[u]/c_1[v,w]$. Since $\X$ is assumed to be
irreducible we have $\overline{\X}\neq 0$, and by Theorem
\ref{rentshcler} there exists $\vartheta\in \K[u]/c_1[v,w]$ such
that $\K[u]/c_1[v,w]^{\overline{\X}}=\K[u]/c_1[\vartheta]$.

\mni Clearly,  $\K[u,p,s]/\ci_{c_1}^{\X}$ is a
$\K[u]/c_1$-subalgebra of $\K[u]/c_1[v,w]$ and we have
$\overline{\X}(p)=\overline{\X}(s)=0$ in $\K[u]/c_1[v,w]$. This
proves that $\K[u,p,s]/\ci_{c_1}^{\X}$ is in fact a
$\K[u]/c_1$-subalgebra of $\K[u]/c_1[\vartheta]$, and as a
consequence there exist polynomials $a(t),b(t)\in \K[u]/c_1[t]$
such that $p=a(\vartheta)$ and $s=b(\vartheta)$ in
$\K[u]/c_1[\vartheta]$. To prove that $a(t)$ is nonconstant we
will prove that $\K[u,p]\cap \ci_{c_1}^{\X}=(c_1)$. Let $k(u,p)\in
\K[u,p]\cap \ci_{c_1}^{\X}$ and write $k(u,p)=c_1(u)p_1(u,v,w)$.
Since $\K[u,p]$ is factorially closed in $\K[u,v,w]$ we have
$p_1(u,v,w)=p_2(u,p)$, and so $k(u,p)=c_1(u)p_2(u,p)$.

Now if $a(t)$ is constant, say $a_0(u)$, then $p-a_0(u)=0$ in
$\K[u]/c_1[v,w]$ and so $p-a_0(u)\in \K[u,p]\cap \ci_{c_1}^{\X}$.
This contradicts the fact that $\K[u,p]\cap \ci=(c_1(u))$.

The fact $\K[u,p]\cap \ci_{c_1}^{\X}=(c_1)$ implies that the
polynomial algebra $\K[u]/c_1[p]$ is a $\K[u]/c_1$-subalgebra of
$\K[u]/c_1[v,w]$. Let us write $a(t)=a_m(u)t^m+\cdots +a_0(u)$
with $m\geq 1$ and $a_m$ a unit of $\K[u]/c_1$. The fact that
$a(\vartheta)-p=0$ in $\K[u]/c_1[v,w]$ implies that $\vartheta$ is
integral over $\K[u]/c_1[p]$. From $s=b(\vartheta)$ in
$\K[u]/c_1[v,w]$ we deduce that $s$ is integral over
$\K[u]/c_1[p]$ as well. Since $\K[u]/c_1[p]$ is UFD and
$\K[u]/c_1[v,w]$ is a domain there exists a unique irreducible
polynomial $h_1(u,p,t)$ which is monic with respect to $t$ such
that $h_1(u,p,s)=0$ in $\K[u]/c_1[v,w]$. Moreover, any other
polynomial $h(u,p,t)$ such that $h(u,p,s)=0$ in $\K[u]/c_1[v,w]$
is a multiple of $h_1$. This means exactly that
$c_1\K[u,v,w]\cap\K[u,p,s]=(c_1,h_1)$ and that $h_1$ is unique, up
to a multiplication by a constant in $\K[u]/c_1$, when viewed as
polynomial in $\K[u]/c_1[p,s]$.

\medskip Now let $a\in
\ci_{c_1}^{\X}$, and notice that in this case reducing $a$ by $h$,
with respect to the lex-order $u\prec p\prec s$, is the same as
performing the Euclidean division of $a$ by $h$ with respect to
$s$. We may thus write $a=qh+r$, with $\deg_s(r)<\deg_s(h)$. Since
$r\in \ci_{c_1}^{\X}$ we may write $r=b_1h+b_2c_1$, and even if it
means reducing $b_2$ by $h$ we may assume that
$\deg_s(b_2)<\deg_s(h)$. By comparing degrees with respect to $s$
in both sides of the last equality we get $r=b_2c_1$, and so $a$
reduces to $0$ by using $h,c_1$. This means exactly that $c_1,h_1$
is a Gr\"obner basis of $\ci_{c_1}^{\X}$ with respect to the
lex-order $u\prec p\prec s$.

\medskip $ii)$ Let $h\in \ci_{c_1}^{\X}$ be a coordinate of
$\K[u]/c_1[p,s]$, and write $h=ac_1+bh_1$. Then over $\K[u]/c_1$
we have $h=bh_1$, and the fact that $h$ is a coordinate of
$\K[u]/c_1[p,s]$ implies in particular that it is irreducible.
This shows that $b$ is a unit of $\K[u]/c_1[p,s]$, and so a
nonzero element of the field $\K[u]/c_1$. As a consequence of
this, $h_1$ is a coordinate of $\K[u]/c_1[p,s]$. The converse is
clear.\end{proof}

\medskip We can now state the main result of this paper.

\begin{theorem}\label{main_theorem}Let $\X$ be a rank
two irreducible locally nilpotent $\K$-derivation of $\K[x,y,z]$,
and let $\K[x,y,z]^{\X}=\K[u,p]$ where $u$ is a coordinate of
$\K[x,y,z]$. Let $s$ be a minimal local slice of $\X$ and write
$\X(s)=c(u)=c_1^{n_1}\cdots c_r^{n_r}$, where the $c_i$'s are
irreducible and pairwise distinct. Then $\X$ is triangulable if
and only if for any $i=1,\ldots,r$ the following hold:

i) the reduced Gr\"obner basis of $\ci_{c_i}^{\X}$ with respect to
the lex-order $u\prec p\prec s$ is $c_i,h_{i}$, where
$h_{i}=Q_i(u,s+\ell_i(u,p))+\mu_i(u)p \mod c_i$ and $\mu_i(u)$ is
a unit $\mod c_i$,

ii) if $\ell(u,p)$ is such that $\ell(u,p)=\ell_i(u,p))\mod c_i$,
and $v=s+\ell(u,p)$ then $u,v$ is a system of coordinates of
$\K[x,y,z]$.

\mni In this case the ideal $\ci_{c}^{\Y}$, where
$\Y=\jac_{(x,y,z)}(u,v,.)$, contains a polynomial of the form
$p+Q(u,s+\ell(u,p))$ and if we let $p+Q(u,s+\ell(u,p))=c(u)w$ then
$u,v,w$ is a coordinate system of $\K[x,y,z]$ which satisfies
$$\X(u)=0,\; \X(v)=c(u),\;\X(w)=\partial_vQ(u,v).$$
\end{theorem}
\begin{proof}``$\Rightarrow$" By Theorem
\ref{triangulability_theorem}, the ideal $\ci_{c}^{\X}$ contains a
polynomial $h^{\star}$ of the form
$p+Q^{\star}(u,s+\ell^{\star}(u,p))$ and if we let
$v^{\star}=s+\ell^{\star}(u,p)$ and $h^{\star}=c(u)w^{\star}$ then
$u,v^{\star},w^{\star}$ is a coordinate system of $\K[x,y,z]$.

\mni For any $i=1,\ldots,r$, let
$h_i^{\star},Q_i^{\star},\ell_i^{\star}$ be respectively the
reductions modulo $c_i$ of $h^{\star},Q^{\star},\ell^{\star}$. The
fact that reduction modulo $c_i$ is a $\K$-algebra homomorphism
implies that $h_i^{\star}=p+Q_i^{\star}(u,s+\ell_i^{\star}(u,p))
\mod c_i$.

\mni Since $h^{\star}$ is a coordinate of $\K[u][p,s]$ it is a
coordinate of $\K[u]/c_1[p,s]$ according to Theorem
\ref{residual_coordinate}. By Lemma \ref{monic_lemma} $i)$ let
$c_i,h_i$ be the reduced Gr\"obner basis of $\ci_{c_i}^{\X}$ with
respect to the lex-order $u\prec p\prec s$. According to Lemma
\ref{monic_lemma} $ii)$ there exists a unit $\nu_i$ modulo $c_i$
such that $h_i^{\star}=\nu_i(u)h_i$. If we let $\mu_i(u)$ be such
that $\mu_i\nu_i=1 \mod c_i$ then we have
$h_i=Q_i(u,s+\ell_i^{\star}(u,p))+\mu_i(u)p$, where $Q_i\in
\K[u,t]$.

\medskip Now let $\ell(u,p)$ be such that
$\ell(u,p)=\ell_i^{\star}(u,p)\mod c_i$ for any $i=1,\ldots,r$.
Since $\ell_i^{\star}(u,p)=\ell^{\star}(u,p)\mod c_i$ we also have
$\ell(u,p)=\ell^{\star}(u,p)\mod c_i$. We claim that
$v=s+\ell(u,p)$ is a $\K[u]$-coordinate of
$\K[u,v^{\star},w^{\star}]$. Indeed, according to Theorem
\ref{residual_coordinate}, it suffices to show that $v$ is a
coordinate of $\K[u]/d(u)[v^{\star},w^{\star}]$ for any
irreducible polynomial $d(u)\in \K[u]$. Depending on $d(u)$ we
have the two following cases.

\mni -- For some $i=1,\ldots,r$, $d(u)$ and $c_i(u)$ are
associate. In this case we have $v=v^{\star}$ in
$\K[u]/d(u)[v^{\star},w^{\star}]$, so $v$ is a coordinate in
$\K[u]/d(u)[v^{\star},w^{\star}]$.

\mni -- For any $i=1,\ldots,r$, $\gcd(d,c_i)=1$. In this case
$c(u)$ is a unit of $\K[u]/d(u)$. Let $\overline{\X}$ be the
$\K[u]/d(u)$-derivation of $\K[u]/d(u)[v^{\star},w^{\star}]$
induced by $\X$. Then $\overline{\X}(c^{-1}v)=1$, which proves
according to Theorem \ref{rentshcler} that $v$ is a coordinate of
$\K[u]/d(u)[v^{\star},w^{\star}]$.

\medskip ``$\Leftarrow$" Assume that $i)$ and $ii)$ hold and
let $\Y=\jac(u,v,.)$. By Lemma \ref{mini_slice_derivation}, $\Y$
is locally nilpotent and we have
$\K[u,v^{\star},w^{\star}]^{\Y}=\K[u,v]$. Moreover, $\Y(p)=-c(u)$
and the fact that $v$ is a coordinate of
$\K[u][v^{\star},w^{\star}]$ implies that $\Y$ has a slice $w$. We
therefore have $\Y(p+c(u)w)=0$, and so $p+c(u)w=Q(u,v)$. The fact
that $u,v,w$ is a coordinate system of $\K[x,y,z]$ follows
immediately from Lemma \ref{wright}, and a direct computation
shows that $\X$ has a triangular form in the coordinate system
$u,v,w$.\end{proof}

\begin{remark}\label{nonuniqueness_remark}
Let $\X$ be a triangular $\K$-derivation and write
$$\X(x)=0,\;\X(y)=c(x),\;\X(z)=q(x,y),$$ and let $p=c(x)z-Q(u,v)$
where $q=\partial_yQ$. From Theorem \ref{main_theorem} $ii)$ we
deduce that any $v=y+d(u)\ell(x,p)$, where $d(u)$ is the maximal
square-free factor of $c(u)$, is a coordinate and gives rise to
another coordinate system $x,v,w$ in which $\X$ has a triangular
form with a different polynomial $Q$. Thus, a triangulable
derivation has many, actually infinitely many, triangular forms.
It is also not clear whether there exists a distinguished form
which could serve as a ``normal form". Nevertheless, it should be
noticed that all the triangular forms and their corresponding
coordinate systems are built out of invariants of $\X$, namely
$\S^{\X}$ and the ideals $\ci_{c_i}^{\X}$ where the $c_i$'s are
the primes factors of $c(u)$.
\end{remark}

\medskip Let us now discuss how to computationally check the
conditions $i)$ and $ii)$ of Theorem \ref{main_theorem}. Assume
that condition $i)$ holds and that we have found a polynomial of
the form $p+Q_i(u,s+\ell_i(u,p))$ in each ideal $\ci_{c_i}^{\X}$.
The computation of $\ell(u,p)$ is then just a matter of Chinese
remaindering. On the other hand, from Lemma
\ref{mini_slice_derivation} we know that
$\Y=\jac_{(x,y,z)}(u,v,.)$, where $v=s+\ell(u,p)$, is locally
nilpotent and $\K[x,y,z]^{\Y}=\K[u,v]$. Thus, $v$ is a coordinate
if and only if $\Y$ has a slice. This may be checked by computing
a minimal local slice starting from the local slice $p$, which
reduces to compute a reduced Gr\"obner basis $G$ of
$c(u)\K[x,y,z]\cap \K[u,v,p]$ with respect to the lex-order
$u\prec v\prec p$. In more explicit terms, $v$ is a coordinate if
and only if the computed Gr\"obner basis is of the form
$c(u),p+Q(u,v)$. Notice that in case $v$ is a coordinate, $G$ also
furnishes a polynomial $w$, with $p+Q(u,v)=c(u)w$, which completes
$u,v$ into a coordinate system and the polynomial $Q$ which is
involved in the triangular form of $\X$.

\mni The condition $i)$ is a matter of functional decomposition of
polynomials, and the fact that we are here dealing with monic
polynomials with respect to $s$ makes it almost trivial.

\begin{lemma}\label{decomposition_lemma}Let $c(u)$ be an
irreducible polynomial of $\K[u,v,w]$, $n$ be a positive integer
and $h\in \K[u,v,w]$ be monic with respect to $w$ and write
$$h=w^d+h_{d-1}(u,v)w^{d-1}+\cdots +h_0(u,v).$$ Then the following
are equivalent:

i) $h=Q(u,w+\ell(u,v))$ in $\K[u]/c^n[v,w]$, with $\ell\in
\K[u]/c^n[v]$ and $Q\in \K[u]/c^n[w]$,

ii) $h(u,v,w-\frac{h_{d-1}}{d})$, viewed in $\K[u]/c^n[v,w]$, is a
polynomial of $\K[u]/c^n[w]$.

\mni In this case, we may choose $\ell= \frac{h_{d-1}}{d}$ and
$Q=h(u,v,w-\ell)$.
\end{lemma}
\begin{proof}$i)\Rightarrow ii)$ Let us write
$Q=w^d+q_{d-1}(u)w^{d-1}+\cdots +q_0(u)$. By expanding
$Q(u,w+\ell(u,v))$ an comparing its coefficients with respect to
$w$ to those of $h$ we get $h_{d-1}(u,v)=d\ell(u,v)+q_{d-1}(u)$.
Therefore, $h(u,w-\frac{h_{d-1}}{d})=Q(u,w-\frac{q_{d-1}(u)}{d})$
and this clearly shows that $h(u,w-\frac{h_{d-1}}{d})\in
\K[u]/c^n[w]$.

\medskip $ii)\Rightarrow i)$ Let us write
$h(u,v,w-\frac{h_{d-1}}{d})=Q(u,w)$. Then
$Q(u,w+\frac{h_{d-1}}{d})=h$ and we have the required
decomposition.\end{proof}

\section{Examples}\label{sec:algo_examples}In this section we  give
two examples to illustrate how our algorithm proceeds. All
derivations are given in a Jacobian form, i.e., as $\jac(f,g,.)$,
since in such a form one can algorithmically check whether the
given derivation is locally nilpotent and if so whether its ring
of constants is generated by $f,g$, see \cite{essen_book}. For
implementation we used the Computer Algebra System \Maple release
$10$.

\begin{example}Consider following example from \cite{daigle96}.
$$\begin{array}{ll} f_1=x,\\
g_1=y+\frac{(xz+y^2)^2}{4},\end{array}$$ and let
$\X=\jac_{(x,y,z)}(f_1,g_1,.)=
\partial_zg_1\partial_y-\partial_yg_1\partial_z$.
The derivation $\X$ is locally nilpotent and its kernel is
$\K[f_1,g_1]$. Our algorithm produces $-x$ as a generator of the
plinth ideal $\S^{\X}$ and $s=-xz-y^2$ as a minimal local slice.
The computation of a Gr\"obner basis of $\ci_x^{\X}$ with respect
to the lex-order $x\prec g_1\prec s$ produces then
$x,(s^2-4g_1)^2+16s$, and the polynomial $(s^2-4p)^2+16s$ cannot
be written in the form $\mu g_1+Q(x,s+\ell(x,g_1))$, where $\mu\in
\K^{\star}$. Therefore, $\X$ is not triangulable.
\end{example}

\begin{example}Consider the following polynomials
$$\begin{array}{ll} f_2=2\,x+y+{z}^{2}-2\,zxy+{x}^{2}{y}^{2},\\
g_2=3\,xy+2\,{x}^{2}-2\,zx+2\,{x}^{2}y+{y}^{2}-yz+x{y}^{2}+{z}^{2}y+{z}^{2
}x-{z}^{3}+3\,{z}^{2}xy-\\\quad
2\,zx{y}^{2}-2\,z{x}^{2}y-3\,z{x}^{2}{y}^{2}+{
x}^{2}{y}^{3}+{x}^{3}{y}^{2}+{x}^{3}{y}^{3}-{z}^{2}+2\,
zxy-{x}^{2}{y}^{2},\end{array}$$ and let
$\Y=\jac_{(x,y,z)}(f_2,g_2,.)$. The derivation $\Y$ is locally
nilpotent and its kernel is $\K[f_2,g_2]$. Moreover, our algorithm
produces $f_2$ as a generator of the plinth ideal $\S^{\Y}$ and
$s=z-xy+1$ as a minimal local slice of $\Y$. The computation of a
Gr\"obner basis of $\ci_{f_2}^{\Y}$ with respect to the lex-order
$f_2\prec g_2\prec s$ produces then $f_2,s^2-2s+g_2+1$. If we let
$u=f_2$ and $v=s-1$ then we get $g_2+v^2=f_2w$, where
$w=-y-x+z-xy$. This gives a coordinate system $u,v,w$ such that
$$\Y(u)=0,\;\Y(v)=u,\Y(w)=2v.$$
\end{example}


\begin{thebibliography}{10}

\bibitem{abhyankar-moh75a}
{S. S.} Abhyankar and {T. T.} Moh, \emph{Embeddings of the line in
the plane},
  J. Reine Angew. Math. \textbf{276} (1975), 148--166.

\bibitem{elkahoui2006a}
{M. A.} Barkatou, H.~{El Houari}, and M.~{El Kahoui},
\emph{Characterization of
  rank two locally nilpotent derivations in dimension three}, Submitted to J.
  Algebra (2006).

\bibitem{bass84a}
H.~Bass, \emph{A non-triangulable action of ${G_a}$ on
$\mathbb{A}^3$}, J. Pure
  Appl. Algebra \textbf{33} (1984), 1--5.

\bibitem{berson-essen2000a}
J.~Berson and {A. van den} Essen, \emph{An algorithm to find a
coordinate's
  mate}, J. of Symbolic Computation \textbf{36} (2003), 835--843.

\bibitem{berson2001a}
J.~Berson, {A. van den} Essen, and S.~Maubach, \emph{Derivations
having
  divergence zero on {$R[X,Y]$}}, Israel J. Math. \textbf{124} (2001),
  115--124.

\bibitem{bhat-duta93a}
{S. M.} Bhatwadekar and {A. K.} Dutta, \emph{On residual variables
and stably
  polynomial algebras}, Comm. Alg. \textbf{21} (1993), no.~2, 635--645.

\bibitem{bhatwadekar97a}
{S. M.} Bhatwadekar and {A. K.} Dutta, \emph{Kernel of locally
nilpotent
  {$R$}-derivations of {$R[X,Y]$}}, Trans. Amer. Math. Soc. \textbf{349}
  (1997), no.~8, 3303--3319.

\bibitem{makar_limanov2005}
A.~Crachiola and L.~{Makar-Limanov}, \emph{An algebraic proof of a
cancellation
  theorem for surfaces}, Preprint available at
  \url{http://www.ihes.fr/IHES/Scientifique/Preprint/prepub.php} (2005).

\bibitem{daigle96}
D.~Daigle, \emph{A necessary and sufficient condition for
triangulability of
  derivations of {$k[X, Y,Z]$}}, J. Pure Appl. Algebra \textbf{113} (1996),
  no.~3, 297--305.

\bibitem{daigle98a}
D.~Daigle and G.~Freudenburg, \emph{Locally nilpotent derivations
over a {UFD}
  and an application to rank two locally nilpotent derivations of {$k[X\sb
  1,\cdots,X\sb n]$}}, J. Algebra \textbf{204} (1998), no.~2, 353--371.

\bibitem{daigle-kaliman2004a}
D.~Daigle and S.~Kaliman, \emph{A note on locally nilpotent
derivations and
  variables of k[{X,Y,Z}]}, preprint available at
  \url{http://aix1.uottawa.ca/~ddaigle/}.

\bibitem{elkahoui2004c}
H.~{El Houari} and M.~{El Kahoui}, \emph{Algorithms for
recognizing coordinates
  in two variables over {UFD's}}, Proceedings of ISSAC'2004, ACM Press, 2004,
  pp.~135--140.

\bibitem{elkahoui2005a}
M.~{El Kahoui}, \emph{{UFD's} with commuting linearly independent
locally
  nilpotent derivations}, J. Algebra \textbf{289} (2005), 446--452.

\bibitem{essen_book}
{A. van den} Essen, \emph{Polynomial automorphisms and the
{J}acobian
  conjecture}, Progress in Mathematics, vol. 190, Birkh\"auser Verlag, Basel,
  2000.

\bibitem{essen2004a}
{A. van den} Essen and {P. van} Rossum, \emph{Coordinates in two
variables over
  a {$\Q$}-algebra}, Trans. Amer. Math. Soc. \textbf{356} (2004), no.~5,
  1691--1703.

\bibitem{walcher97a}
{D. R.} Finston and S.~Walcher, \emph{Centralizers of locally
nilpotent
  derivations}, J. Pure Appl. Algebra \textbf{120} (1997), no.~1, 39--49.
  \MR{MR1466096 (98d:14057)}

\bibitem{freudenburg95a}
G.~Freudenburg, \emph{Triangulability criteria for additive group
actions on
  affine space}, J. Pure and Applied Algebra \textbf{105} (1995), 267--275.

\bibitem{freudenburg_book}
G.~Freudenburg, \emph{{Algebraic theory of locally nilpotent
derivations.
  Invariant theory and algebraic transformation groups VII.}}, {Encyclopaedia
  of Mathematical Sciences 136. 260~p.}, 2006 (English).

\bibitem{gathen90}
{J. von zur} Gathen, \emph{Functional decomposition of
polynomials: the tame
  case}, J. Symbolic Comput. \textbf{9} (1990), no.~3, 281--299.

\bibitem{jung42}
{H. W. E.} Jung, \emph{\"{U}ber ganze birationale
{T}ransformationen der
  {E}bene}, J. Reine Angew. Math. \textbf{184} (1942), 161--174.

\bibitem{kaliman2002a}
S.~Kaliman, \emph{Polynomials with general $\mathbb{C}^{2}$-fibers
are
  variables}, Pacific J. Math. \textbf{203} (2002), no.~1, 161--190.

\bibitem{kambayashi75a}
T.~Kambayashi, \emph{On the absence of nontrivial separable forms
of the affine
  plane}, J. Algebra \textbf{35} (1975), no.~449-456.

\bibitem{kulk53}
{W. van der} Kulk, \emph{On polynomial rings in two variables},
Nieuw Arch.
  Wiskunde (3) \textbf{1} (1953), 33--41.

\bibitem{shpilrain2005a}
{C. M.} Lam, V.~Shpilrain, and {J-T.} Yu, \emph{Recognizing and
parametrizing
  curves isomorphic to a line}, To appear in J. Symbolic Computation (2007).

\bibitem{Miyanishi85a}
M.~Miyanishi, \emph{Normal affine subalgebras of a polynomial
ring}, Algebra
  and Topological Theories--to the memory of Dr. Takehiko MIYATA, Kinokuniya
  (1985), 37--51.

\bibitem{popov87}
{V. L.} Popov, \emph{On actions of {${\bf G}\sb a$} on {${\bf
A}\sp n$}},
  Algebraic groups Utrecht 1986, Lecture Notes in Math., vol. 1271, Springer,
  Berlin, 1987, pp.~237--242.

\bibitem{rentschler68}
R.~Rentschler, \emph{Op\'erations du groupe additif sur le plan
affine}, C. R.
  Acad. Sci. Paris S\'er. A-B \textbf{267} (1968), 384--387.

\bibitem{shestakov2004a}
{I. P.} Shestakov and {U. U.} Umirbaev, \emph{The tame and the
wild
  automorphisms of polynomial rings in three variables}, J. Amer. Math. Soc.
  \textbf{17} (2004), no.~1, 197--227.

\bibitem{shpilrain97a}
V.~Shpilrain and J-T. Yu, \emph{Polynomial automorphisms and
{G}r{\"o}bner
  reductions}, J. Algebra \textbf{197} (1997), 546--558.

\bibitem{wright81a}
D.~L. Wright, \emph{On the {J}acobian {C}onjecture}, Illinois J.
Math
  \textbf{25} (1981), no.~3, 423--440.

\end{thebibliography}

\providecommand{\bysame}{\leavevmode\hbox
to3em{\hrulefill}\thinspace}
\providecommand{\MR}{\relax\ifhmode\unskip\space\fi MR }
\providecommand{\MRhref}[2]{%
  \href{http://www.ams.org/mathscinet-getitem?mr=#1}{#2}
} \providecommand{\href}[2]{#2}

\end{document}